\documentclass[conference]{IEEEtran}
\IEEEoverridecommandlockouts
\usepackage{cite}
\usepackage{amsmath,amssymb,amsfonts}
\usepackage{algorithmic}
\usepackage{graphicx}
\usepackage{textcomp}
\usepackage{xcolor}
\usepackage{makecell}
\def\BibTeX{{\rm B\kern-.05em{\sc i\kern-.025em b}\kern-.08em
    T\kern-.1667em\lower.7ex\hbox{E}\kern-.125emX}}
\begin{document}

\title{Data-driven Optimization Approaches to the Power System Planning under Uncertainty *
{\footnotesize }

}

\author{\IEEEauthorblockN{Shuhan Lyu}
\IEEEauthorblockA{\textit{Department of Electrical Engineering and Computer Science} \\
\textit{Technical University of Berlin}\\
Berlin, Germany\\
shuhan.lyu@campus.tu-berlin.de}
}

\maketitle

\begin{abstract}
In order to protect the environment and address fossil fuel scarcity, renewable energy is increasingly used for power generation. However, due to the uncertainties it brings to electricity production, deterministic optimization is no longer sufficient for operational needs. Therefore, a large number of optimization techniques under uncertainty have been proposed, which provide good ways to address uncertainties. This paper selects three of the more important optimization techniques under uncertainty to introduce: stochastic programming (SP), robust optimization (RO), and a novel approach named distributionally robust optimization (DRO) based on the first two. We explain each approach’s basic framework and general process using specific examples. The focus is on how each method addresses the uncertainties. In addition, we also compare their strengths and weaknesses and discuss future research directions.
\end{abstract}


\section{Introduction}
Since electricity was invented, it has occupied an increasingly important place in our lives. Today, we can no longer imagine living without electricity. The traditional way of generating electricity is thermal power, which relies on burning coal to produce electricity. However, due to the shortage of fossil fuels and the huge pollution it brings to the environment, such as a large amount of carbon dioxide emissions, people are committed to using more renewable energy for electricity generation today. For example, wind energy is one of the most commonly used green resources on the market today. 

In addition to energy updates, we also have some novel power generation systems, such as the Microgrid. The Berkeley Lab defines it as: “A microgrid consists of energy generation and energy storage that can power a building, campus, or community when not connected to the electric grid, e.g. in the event of a disaster \cite{b1}."  Compared with the bulk power system, a microgrid is more flexible and can be used to supply power in emergencies and remote areas. Besides, a microgrid can effectively integrate various sources of distributed generation, especially renewable energy sources \cite{b2}.  Therefore, it can offer an option to balance the need to reduce carbon emissions with continuing to provide reliable electric energy in periods when renewable power sources are unavailable \cite{b4} \cite{b5}.

Nevertheless, every coin has two sides. While renewable energy has many advantages over traditional energy, it also brings many uncertainties, which makes it difficult to optimize the short-term operation of power systems\cite{c1}\cite{c2}.
In particular, the inherent uncertainty of wind energy makes it very difficult to be predict. The output is fluctuating and we should also take into account occasional weather extremes, such as hurricanes \cite{b6} and sandstorms. On the other hand, the consumption demand also changes over time, for microgrids when the demand is greater than the net load, we must purchase additional power from the electricity market to fill the demand \cite{b7}. The ensuing question is, the electricity market price also has volatility, which means besides the power generation cost, we should also consider the purchase cost. Moreover, given the restricted scale and size of a microgrid, even minor errors in the wind and solar power projections can cause significant uncertainty in real-time operations \cite{b3}. Therefore, it is particularly important to develop a suitable power generation plan that meets demand while minimizing expenses.

The proper term for this optimal problem is “unit commitment” (UC). The UC problem refers to the mathematical problem of scheduling the commitment, start-up, and shutdown statuses of dispatchable generators, such as thermal resources. Generally it falls into two categories: deterministic optimization, and  optimization techniques under uncertainty. If it were possible to know all of the parameters, deterministic optimization models could be used for solving UC problems, but if there is uncertainty, we should use optimization techniques under uncertainty. Classification by time, we have long-term and short-term planning, long-term means week-ahead or even year-ahead, while short-term means hours-ahead to day-ahead. This paper is to review the works that have contributed to the modeling for short-term UC focusing on uncertainties of renewables. We introduce three main modeling methods within this classification, namely stochastic programming (SP), robust optimization (RO), and distributionally robust optimization (DRO). 

The remainder of this paper is organized as follows: In Section \uppercase\expandafter{\romannumeral2} we first explain  optimization techniques under uncertainty in general, then introduce the two basic models in this category, namely stochastic programming and robust optimization, in turn. We present the novel distributionally robust optimization model based on the first two in Section \uppercase\expandafter{\romannumeral3} and discuss the future research directions in Section \uppercase\expandafter{\romannumeral4}. Conclusions are drawn in Section \uppercase\expandafter{\romannumeral5}.

\section{Two Basic UC Models}
While deterministic UC optimization assumes the next-day situation is known, stochastic models consider the uncertainties for the following day or hours \cite{b9}. In most cases the model for  optimization techniques under uncertainty is two-stage. In the first stage, we make the day-ahead decision, namely when and which generating units at each power station are startup or shutdown. Because the traditional units like coal generators cannot be turned on or off quickly in real time, we must make the decision ahead of time. Then in the second stage, we determine the real-time plan. It means to decide the individual power outputs of the scheduled generating units at each time period involving uncertainty created by renewable energy. Hence, one of the most important points in  optimization techniques under uncertainty is to deal with uncertainties. Next, we introduce the two basic models, stochastic programming and robust optimization, in turn. We first explain their basic models and processes, then compare the advantages and disadvantages of the two methods to have a more comprehensive understanding.

\subsection{Stochastic Programming}

The basic concept of stochastic programming (SP) is scenario-based uncertainty representation. Each scenario depicts a possible realization, according to which using more scenarios means getting higher quality results. 

The general process is divided into three steps:

1. build a probability distribution $P_0$ based on historical data

2. use Monte Carlo simulation for populating scenarios subject to $P_0$

3. calculate the optimal solution with existing scenarios

In the first step, by analyzing the historical data, we can assume that wind power follows a specific distribution, such as the normal distribution \cite{b6} or the Weibull distribution \cite{b7}. The wind power prediction is formulated with the expected value and the variance from $P_0$. Then we generate more scenarios according to $P_0$ for higher quality results. Each scenario is assigned a random wind power output based on the wind power prediction.

The basic model can be formulated as follows \cite{b9}:

\begin{equation}
  \min_{u\in\mathbb{U}} c^\tau u + E_\xi [F(u,\xi)]  \label{4}
\end{equation}
\begin{equation}
    F(u,s) = \min_{p_s,f_s}  f(p_s) \label{5}
\end{equation}

$u$ refers to commitment decisions of tradition units, such as coal generators. $\mathbb{U}$ consists of all possible commitment decisions and the vector $c$ stores various startup and shutdown expenses. As previously stated, due to the nature of traditional units, $u$ should be decided ahead of time. That means in the first-stage ($i.e.$,  \eqref{4}), we make the day-ahead decision $u$.

The second term in \eqref{4} is the expected real-time operation cost, and $\xi$ is the random vector that has a known joint probability distribution ($P_0$), which is made based on the scenarios before. For each realization $s$ of the random vector $\xi$, the second-stage problem can be formulated as \eqref{5}. $p_s$ comprises dispatches as well as reserves in multiple periods, and $f_s$ is the vector of other second-stage decisions (such as power flows and pump-storage hydro unit decisions). The goal of the second-stage problem is to find the minimum of the fuel cost function $f(p_s)$.

 To compute the solution, we use decomposition methods since, in the second-stage problem, different scenarios (different realizations $s$ of $\xi$) are not inextricably bound to each other. Therefore, we can handle them separately, which means grouping them into smaller independent optimization problems \cite{b6}. Commonly used approaches are Benders decomposition and L-shaped algorithm \cite{b10}.

Stochastic programming is always criticized for not taking into account the lack of power supply that may be caused by unexpected situations but in fact, this issue can be solved by applying some risk measures. For instance, \cite{b11} uses Expected Load Not Served (ELNS), and \cite{b12} utilizes Loss of Load Probability (LOLP). 

The ELNS value is evaluated at each time period by taking the expectation of the total net load minus the total dispatch.
The main advantage of ELNS is its easiness of calculation and hence it can be included in the objective function or a bounding constraint. However, since ELNS is based on expected values, it cannot tell how risky a certain decision may be. LOLP can be a good solution to this problem, which is defined as follows \cite{b9}:
\begin{equation}
LOLP_t = Prob(1^\tau p_t^\xi\geq d^\xi_t)\label{1}
\end{equation}
where $1^\tau p_t^\xi$ is the overall dispatch susceptible to sporadic failures and $d^\xi_t$ is the random demand in period $t$. 
It calculates the probability that the total dispatch subject to random outages or other uncertainties meets the actual demand. To ensure system reliability, we can restrict the LOLP value to be greater than a confidence level $1 - \varepsilon$ for all time periods $t$ and $\varepsilon$ should be a very small number, then add the constraints to the model. For instance, the constraint can be described as follows \cite{b9}:

\begin{equation}
    Prob (1^\tau p_t^\xi + l_t^\tau\geq d^\xi_t) \geq 1-\varepsilon, \forall t
\end{equation}

where $l_t^\tau$ represents the load shedding necessary to achieve the power balance considering uncertainties. The probabilistic constraints can be utilized to restrict the total load loss or load losses at different zones \cite{b12}. 

The main advantage of LOLP is the use of the explicit reliability measure, with which we are able to balance total cost and reliability. However, it also comes with significant computational challenges. For example, uncertain outages could mean a large number of scenarios, for each of which a binary variable needs to be introduced, with 0 and 1 representing if a scenario will be selected or not. Thus, a lot of power outages mean a huge amount of computing. As a result, approximations based on a limited number of outages are used \cite{b13}.

\subsection{Robust Optimization}
Instead of assuming the probability distribution, the robust optimization represents uncertainty with an uncertainty set, which defines the upper and lower bounds of the uncertainty. It is also worth mentioning that robustness means considering the worst case. 

The general process can be described in three steps:

1. construct an uncertainty set $S$ (usually a set of percentiles) using historical data

2. define an interval for each time period with lower and upper bounds

3. calculate the solution in consideration of the worst case

There are many different ways to build the uncertainty set and to define the worst case, we introduce one of them which is presented in \cite{b14}.

Suppose we need to define an interval $[Q_l, Q_u]$ of wind power $q_1$ in time period $t_1$. As previously stated, the uncertainty set is always a set of percentiles. For instance, we can choose the 0.05- and 0.95-percentile from the uncertainty set $S$ as the lower and upper bounds for $q_1$. As a result, $Q_l$ is equal to the 0.05-percentile and $Q_u$ the 0.95-percentile.
When $q_1$ is greater than the 0.95-percentile, we should curtail the excess wind power, and if the $q_1$ is smaller than the 0.05-percentile, we need to curtail the partial loads. We do this for every time period and the remaining question is, how do we define the worst case? If we define it as where the wind output is equal to the boundary value, which means $q_1 = Q_l$ or $q_1 = Q_u$, 
we would be too conservative. It is better to formalize it as where the wind output has a huge fluctuation. To control the conservativeness, we use an integer $T$ for the number of time periods that allow $q_t$ to be far away from the mean value of the defined interval. For example, if $T$ is equal to 0, every $q_t$ should be close to the mean value, and if $T$ is equal to 12, $q_t$ could take any value within the interval in 12 periods. Hence, a larger $T$ means that we have to take into account more worst cases, implying a more conservative result.

The general model of robust optimization can be defined as \cite{b9}:
\begin{equation}
    \min_{u\in\mathbb{U}} \left\{ c^\tau u + \max_{v\in\mathbb{V}} [F(u,v)] \right\} \label{10}
\end{equation}

where $u$, $\mathbb{U}$, and the vector $c$ have the same definitions as in the stochastic programming model. In the first stage, we also do the day-ahead decision $u$. 
$\mathbb{V}$ is the predefined uncertainty set and $v$ the deterministic parameter. The function $F(u,v)$ speaks of the real-time dispatch cost. It is worth mentioning that $v$ is not a parameter vector but a variable vector, because the purpose of \eqref{10} is to obtain the minimum total cost for the worst case. $F(u,v)$ is defined as the optimal value for the minimization problem shown in \eqref{11}.

\begin{equation}
    F(u,v) = \min_{p,f} q^\tau p \label{11}
\end{equation}

$q$ is the coefficient vector, as the piecewise linear approximation is used. The definitions of $p$ and $f$ are the same as in \eqref{5}.

The RO problem is a min-max-min problem, which can not be directly solved using existing software. Thus, we need to develop novel solution frameworks for it. Separate economic dispatch constraints corresponding to each outcome must be taken into account. Therefore, we need to use decomposition methods, such as Benders decomposition, to split the problem into multiple sub-problems. Moreover, since the RO approach takes into account the worst-case scenario, iterations should also be used to solve the problem.

\subsection{Comparison}
In summary, comparing stochastic programming and robust optimization, each method has its advantages and disadvantages, which are listed in Table \ref{tab1}.

For stochastic programming, the predefined probability distribution could be unreliable which leads to an inaccurate result. Although a huge number of scenarios always results in a prohibitively high computational burden \cite{b9}, a significant strength of SP is that there is no need to create a new framework to compute the result. Because the SP model is not a min-max-min model, the problem can be directly solved with existing methods and software, such as CPLEX.

In contrast, we need different expertise and rationale for different uncertainty sets. To compute the solution, we always need to develop novel frameworks using new algorithms, but the computational demand is correspondingly less. Besides, since we do not need to define a static $P_0$, the results of robust optimization are more reliable compared to stochastic programming. However, because RO always produces over-conservative solutions, the expected total expenses could be very high.

\begin{table}[htbp]
\caption{Comparision Of Methods}
\begin{center}
\begin{tabular}{|c|c|c|}
\hline
\textbf{Methods}&\textbf{Advantages} & \textbf{Disadvantages} \\
\hline
\makecell{Stochastic\\Programming} & \makecell[l]{1. Easier to understand than\\minimizing the worst-case\\cost\\2. May use risk measures to\\solve robustness concerns}& \makecell[l]{1. Static probability\\distribution assumption\\2. High computational\\ burden for a large\\ number of scenarios}\\
\hline
\makecell{Robust\\Optimization} & \makecell[l]{1. No assumption of the\\probability distribution\\2. Reduce the computational\\load} & \makecell[l]{1. Over-conservative\\solutions\\2. Need expertise and\\ rationale for uncertainty\\sets}\\
\hline

\end{tabular}
\label{tab1}
\end{center}
\end{table}

\section{distributionally robust optimization}
Considering the disadvantages of both methods, there is a novel UC modeling approach that is more reliable than SP and not as conservative as RO, namely the distributionally robust optimization, “robust” means it also takes into consideration the worst case. In this method, an ambiguity set of probability distributions is used to represent the uncertainty. The framework can be described as follows \cite{b15}: 

\begin{equation}
    \min_{u\in\mathbb{U}} c^\tau u + \sup_{\mathbb{P}\in\mathcal{P}} E_\mathbb{P}[f(u,\omega)] \label{6}
\end{equation}

where $u$, $\mathbb{U}$, and $c$ are the same defined as in \eqref{4} and in \eqref{10}. 
In the first stage, we make the decision $u$ for the traditional generators. 
$\omega$ is a random variable following the probability distribution $\mathbb{P}$.  
The ambiguity set $\mathcal{P}$ includes other probability distributions around the nominal probability distribution based on a certain ambiguity set construction method and a specified set size, in this case it contains all possible probability distributions of $\omega$. 

The objective function $F(u,\omega)$ aims to find the minimum expected value under the worst-case distribution over $\mathcal{P}$.

\begin{equation}
    F(u,\omega) = \max_{m\leq M} F_m (u,\omega) \label{7}
\end{equation}

We can see from \eqref{7} that $F(u,\omega)$ is defined as the point-wise maximum of $F_m$, which represents the traditional unit dispatch cost. 

The DRO model is always a “min-max-min” programming model, the goal of the first phase is to get the minimum total cost, and in the second stage, we aim to obtain the optimal solution for the dispatch problem. That means the maximum usage with minimum operation cost, which is the last “min” refers to.

Generally, the process has three steps:

1. form a nominal probability distribution $P_0$ from historical data

2. build the ambiguity set $\mathcal{P}$ based on $P_0$

3. solve the problem considering the worst-case scenario

The literature documents a number of methods for constructing $P_0$ in the first step. 
In \cite{b133} and \cite{b144}, they first apply k-means clustering to group the days that contain the required information in order to construct the nominal probability distribution. In \cite{b155}, they introduce both parametric and non-parametric methods for $P_0$, while in \cite{b15}, $P_0$ is obtained with the Dirac measure.

In the second step, there are also many approaches for building the ambiguity set, we present two popular ways here: KL divergence and Wasserstein metric. 
The models of KL divergence can be formulated as follows \cite{b155}:

\begin{equation}
   \mathbb{P} = \left\{ P\in\mathbb{D}:D_{KL}(P||P_0)\leq \eta\right\} \label{2}
\end{equation}

$\eta$ is the divergence tolerance, $P_0$ is the nominal distribution, and $\mathbb{D}$ represents the set of all probability distributions. $D_{KL}(P||P_0)$ is the KL divergence of $P$ from $P_0$. To build the ambiguity set, we collect all $P$ within $\mathbb{D}$, whose KL divergence from $P_0$ is less than or equal to $\eta$. In other words,  $\eta$ control the conservativeness of the ambiguity set.
If the tolerance parameter $\eta$ equals zero, then the ambiguity set includes only the nominal distribution $P_0$, and if it is close to infinity, the ambiguity set contains all probability distributions. 

To determine the divergence tolerance, we could use, for example, the confidence level. Then the model could be described as \cite{b155} :

\begin{equation}
    \mathbb{P} = \left\{ P\in\mathbb{D}:D_{KL}(P||P_0)\leq \eta_\alpha \right\} \label{2}
\end{equation}

where $\alpha$ is the confidence level and $\eta_\alpha$ is the divergence tolerance with $\alpha$. The constraint describing this ambiguity set $\mathbb{P}$ guarantees that the true distribution is within $\mathbb{P}$ with $\alpha$ confidence level \cite{b155}.

The following is the model of the Wasserstein metric \cite{b15}:

\begin{equation}
\mathcal{P}:=\left\{\mathbb{P}|W(\mathbb{P}_N,\mathbb{P})\leq\epsilon N\right\}\label{3}
\end{equation}

$\mathcal{P}$ represents a Wasserstein ball with the radius $\epsilon N$ and the center $\mathbb{P}_N$. $W(\mathbb{P}_N,\mathbb{P})$ is defined as the Wasserstein distance between the marginal distribution $\mathbb{P}_N$ and $\mathbb{P}$. Similar to the KL divergence, in the Wasserstein metric, we also select all $\mathbb{P}$, which are close enough to $\mathbb{P}_N$. The radius $\epsilon N$ is also used to control the conservativeness of the model, which means to determine the size of the ambiguity set. We can  determine the value of the tolerance parameter in different manners, beside the confidence interval, there are other optimal solutions for the tolerance parameter. System operators can decide  which method to use according to their needs.

Putting the two methods together, we could see that the basic ideas of both methods are generally consistent. They first measure the divergence/distance between $P$ and $P_0$ in their ways, then set a specific value for the conservativeness. When it comes to the biggest difference between the two methods, the KL-Divergence-based ambiguity set can only be produced from data when the underlying distributions are supported on a finite range of values. Nevertheless, in fact, the underlying distributions of wind power have continuous support\cite{b15}. The computation of the Wasserstein- metric-based technique, on the other hand, is proportionally more difficult. The exact method used to construct the ambiguity depends on the situations.

 In \cite{b133}-\cite{b155} they use KL divergence to solve the UC problem with different uncertainty resources. In \cite{b133} the uncertainty is the net load, in \cite{b144} the microgrids, and in \cite{b155} the wind power. \cite{b15} presents a Wasserstein metric-based DRO under wind uncertainty.

Since the DRO model is a three-level model that cannot be solved directly with existing frameworks, we always reformulate it into a simpler model, for instance a one-level model, and apply the Benders decomposition or even two-level decomposition to reduce the difficulty of calculation. Alternatively, we could also use approximate computation and then confirm the condition under which the approximation model becomes precise \cite{b15}.

\section{Future Research Direction}
Future research should focus more on providing accurate weather forecasts, especially renewable generation forecasting. How to incorporate it into system operation to deliver more optimized results is also an important question. Moreover, it is very much the key component to overcoming the computational challenge. To reduce the computational burden, we have different coping methods for different models. 
 For the stochastic programming model, the effective methods are efficient scenario selection and reduction, which can remove scenarios without affecting the quality of the results. When it comes to robust optimization, we need advanced algorithms for more general uncertainty sets. Besides, existing unit commitment software packages, such as \cite{c3}, can be extended to incorporate uncertainty.The new paradigm of contextual stochastic optimization could also be applied to take unit commitment decisions under uncertainty, as proposed in \cite{c4}.
  Furthermore, it is important to develop more new frameworks, which can simplify the calculating process and give efficient solutions for distributionally robust optimization modeling. Last but not least, multi-thinking is also a critical point for this optimization problem, because uncertainties stem from many sources, from changes in electricity demand to fluctuations in renewable energy output. Today, sustainable development gains more and more attention, and new policies such as carbon tax should also be considered and added to the constraints \cite{b166}.

\section{Conclusion}
In order to reduce the environmental pollution caused by thermal power generation, to address fossil fuel scarcity, as well as to meet the growing electricity demand, renewable energy is increasingly used in the market today. To overcome the instability of renewable energy sources, scientists have worked on more flexible modeling methods. In this paper, we select three of the main methods to introduce: stochastic programming, robust optimization, and distributionally robust optimization. The novel DRO is relatively the most effective, which overcomes the shortcomings of the first two methods and takes the best from the worst. It is more reliable than SP while providing results that are not as conservative as RO. However, in actual production, the exact method to be used still depends on the calculation requirements that can be met and the decision maker's need for conservativeness. Each method has its most applicable case.


\begin{thebibliography}{00}
\bibitem{b1} Microgrids and Vehicle-Grid Integration. Berkeley Lab. Retrieved 21 June 2022.
\bibitem{b2} J. Hu and A. Lanzon, "Distributed Finite-Time Consensus Control for Heterogeneous Battery Energy Storage Systems in Droop-Controlled Microgrids," in IEEE Transactions on Smart Grid, vol. 10, no. 5, pp. 4751-4761, Sept. 2019.

\bibitem{b4} M. Saleh, Y. Esa and A. A. Mohamed, "Communication-Based Control for DC Microgrids," in IEEE Transactions on Smart Grid, vol. 10, no. 2, pp. 2180-2195, March 2019. 
\bibitem{b5} D. E. Olivares $et$ $al$., "Trends in Microgrid Control," in IEEE Transactions on Smart Grid, vol. 5, no. 4, pp. 1905-1919, July 2014.

\bibitem{c1} O. Yurdakul, F. Eser, F. Sivrikaya and S. Albayrak, "Very Short-Term Power System Frequency Forecasting," in IEEE Access, vol. 8, pp. 141234-141245, 2020.

\bibitem{c2}O. Yurdakul, A. Meyer, F. Sivrikaya and S. Albayrak, "Forecasting Daily Primary Three-Hour Net Load Ramps in the CAISO System," 2020 52nd North American Power Symposium (NAPS), 2021, pp. 1-6.


\bibitem{b6} J. Wang, M. Shahidehpour and Z. Li, "Security-constrained unit commitment with volatile wind power generation," 2009 IEEE Power and Energy Society General Meeting, 2009, pp. 1-1.
\bibitem{b7}B. Venkatesh, P. Yu, H. Gooi, and D. Choling, “Fuzzy MILP unit commitment incorporating wind generators,” IEEE Trans. Power Syst., vol.23, no. 4, pp. 1738–1746, Nov. 2008.
\bibitem{b3}W. Su, J. Wang and J. Roh, "Stochastic Energy Scheduling in Microgrids With Intermittent Renewable Energy Resources," in IEEE Transactions on Smart Grid, vol. 5, no. 4, pp. 1876-1883, July 2014.
\bibitem{b9} Q. P. Zheng, J. Wang and A. L. Liu, "Stochastic Optimization for Unit Commitment—A Review," in IEEE Transactions on Power Systems, vol. 30, no. 4, pp. 1913-1924, July 2015.
\bibitem{b10}  J. R. Birge and F. Louveaux, Introduction to Stochastic Programming. New York, NY, USA: Springer Verlag, 1997.
\bibitem{b11}H. B. Gooi, D. P. Mendes, K. R. W. Bell and D. S. Kirschen, "Optimal scheduling of spinning reserve," in IEEE Transactions on Power Systems, vol. 14, no. 4, pp. 1485-1492, Nov. 1999.
\bibitem{b12}A. Ahmadi-Khatir, M. Bozorg and R. Cherkaoui, "Probabilistic Spinning Reserve Provision Model in Multi-Control Zone Power System," in IEEE Transactions on Power Systems, vol. 28, no. 3, pp. 2819-2829, Aug. 2013.
\bibitem{b13}F. Bouffard and F. Galiana, “An electricity market with a probabilistic
spinning reserve criterion,” IEEE Trans. Power Syst., vol. 19, no. 1, pp.300–307, Feb. 2004.
\bibitem{b14}R. Jiang, J. Wang and Y. Guan, "Robust Unit Commitment With Wind Power and Pumped Storage Hydro," in IEEE Transactions on Power Systems, vol. 27, no. 2, pp. 800-810, May 2012.


\bibitem{b15} R. Zhu, H. Wei and X. Bai, "Wasserstein Metric Based Distributionally Robust Approximate Framework for Unit Commitment," in IEEE Transactions on Power Systems, vol. 34, no. 4, pp. 2991-3001, July 2019.

\bibitem{b133}O.Yurdakul, F.Sivrikaya and S.Albayrak, "Kullback-Leibler Divergence-Based Distributionally Robust Unit Commitment Under Net Load Uncertainty," 2021 IEEE Madrid PowerTech, 2021, pp. 1-6.
\bibitem{b144} O. Yurdakul, F. Sivrikaya and S. Albayrak, "A Distributionally Robust Optimization Approach to Unit Commitment in Microgrids," 2021 IEEE Power and Energy Society General Meeting (PESGM), 2021, pp. 1-5.
\bibitem{b155}Y. Chen, Q. Guo, H. Sun, Z. Li, W. Wu and Z. Li, "A Distributionally Robust Optimization Model for Unit Commitment Based on Kullback–Leibler Divergence," in IEEE Transactions on Power Systems, vol. 33, no. 5, pp. 5147-5160, Sept. 2018.

\bibitem{c3} Alinson S. Xavier, Aleksandr M. Kazachkov, Ogün Yurdakul, Feng Qiu, "UnitCommitment.jl: A Julia/JuMP Optimization Package for Security-Constrained Unit Commitment (Version 0.3)," Zenodo (2022), 2022.

\bibitem{c4}  O. Yurdakul, F. Qiu, and S. Albayrak,  "Predictive Prescription of Unit Commitment Decisions Under Net Load Uncertainty. " arXiv, 2022.

\bibitem{b166} O. Yurdakul, F. Sivrikaya and S. Albayrak, "Quantification of the Impact of GHG Emissions on Unit Commitment in Microgrids," 2020 IEEE PES Transmission \& Distribution Conference and Exhibition - Latin America (T\&D LA), 2020, pp. 1-6
\end{thebibliography}
\end{document}